\documentclass[12pt]{amsart}

\usepackage[dvips]{graphicx}
\usepackage{amsfonts,amssymb,amscd,bm,amsthm}
\newcommand{\Z}{{\mathbb Z}}

\newcommand{\e}{{\varepsilon}}

\newcommand{\al}{{\alpha}}
\newcommand{\be}{{\beta}}
\newcommand{\de}{{\delta}}

\newcommand{\ga}{{\gamma}}

\newcommand{\De}{{\Delta}}

\newcommand{\ti}{\tilde}

\newcommand{\id}{{\mathrm{id}}}

\newcommand{\UD}{{\mathrm{\bf UD}}}
\newcommand{\OD}{{\mathrm{\bf OD}}}
\newcommand{\UC}{{\mathrm{\bf UC}}}
\newcommand{\OC}{{\mathrm{\bf OC}}}
\newcommand{\ub}{{\mathrm{\bf B}}}
\newcommand{\pb}{{\mathrm{\bf P}}}
\newcommand{\nb}{{\mathrm{Nbhd}}}
\newcommand{\supp}{{\mathrm{supp}}}

\newcommand{\bx}{{ \mathbf x }}

\newcommand{\ba}{{ \bm{\alpha} }}
\newcommand{\bb}{{ \bm{\beta} }}
\newcommand{\bg}{{ \bm{\gamma} }}

\newcommand{\bla}{{ \bm{\lambda} }}
\newcommand{\bmu}{{ \bm{\mu} }}
\newcommand{\br}{{ \bm{\rho} }}

\newcommand{\vl}{{\; | \;}}
\newcommand{\es}{{\varnothing}}

\newcommand{\ab}[1]{\langle #1 \rangle}

\newtheorem{theorem}{Theorem}[section]
\newtheorem{corollary}[theorem]{Corollary}
\newtheorem{proposition}[theorem]{Proposition}
\newtheorem{lemma}[theorem]{Lemma}

\theoremstyle{definition}
\newtheorem{definition}[theorem]{Definition}
\newtheorem{example}[theorem]{Example}
\newtheorem{algorithm}[theorem]{Algorithm}

\newtheorem{problem}[theorem]{Problem}

\title[Braid groups of graphs with applications to motion planning]
{Computing braid groups of graphs with \\
 applications to robot motion planning}

\author[kurlin]{V.~Kurlin}
\address{Department of Mathematical Sciences,   
 The University of Durham, \\
 Durham DH1 3LE, United Kingdom}
\email{ vitaliy.kurlin@durham.ac.uk }

\date{First version: 7 August 2009.}

\begin{document}

\begin{abstract}
We design an algorithm writing down presentations of graph braid groups.
Generators are represented in terms of actual motions of 
 robots moving without collisions on a given graph.
A key ingredient is a new motion planning algorithm whose complexity
 is linear in the number of edges and quadratic in the number of robots.
The computing algorithm implies that 2-point braid groups of all light planar graphs 
 have presentations where all relators are commutators.
\end{abstract}

\maketitle


\section{Introduction}
\label{sect:Intro}

\subsection{Brief summary}
\label{subs:Summary}
\noindent
\smallskip

This is a research on the interface between topology and graph theory
 with applications to motion planning algorithms in robotics.
We consider moving objects as zero-size points travelling without collisions 
 along fixed tracks forming a graph, say on a factory floor or road map.
We prefer to call these objects `robots', although the reader 
 may use a more neutral and abstract word like `token'.  
\smallskip

For practical reasons we study discrete analogues of configuration spaces of graphs,
 where robots can not be very close to each other, roughly one edge apart.
This discrete approach reduces the motion planning of real (not zero-size) vehicles
 to combinatorial questions about ideal robots moving on a subdivided graph.


\subsection{Graphs and theirs configuration spaces}
\label{subs:ConfSpaces}
\noindent
\smallskip

First we recall basic notions.
A \emph{graph} $G$ is a 1-dimensional finite CW complex, 
 whose 1-cells are supposed to be open.
The 0-cells and open 1-cells are called \emph{vertices} and 
 \emph{edges}, respectively.
If the endpoints of an edge $e$ are the same then $e$ is called a \emph{loop}. 
A \emph{multiple} edge is a collection of edges with the same distinct endpoints.
The topological \emph{closure} $\bar e$ of an edge $e$ 
 is the edge $e$ itself with its endpoints.
\smallskip

The \emph{degree} $\deg v$ of a vertex $v$
 is the number of edges attached to $v$, 
 i.e. a loop contributes 2 to the degree of its vertex.
Vertices of degrees~1 and~2 are \emph{hanging} and \emph{trivial}, respectively.
Vertices of degree at least 3 are \emph{essential}.
A \emph{path} (a \emph{cycle}, respectively) of length $k$ in $G$ is 
 a subgraph consisting of $k$ edges and homeomorphic to a segment
 (a circle, respectively).
A \emph{tree} is a connected graph without cycles.
\smallskip

The direct product $G^n=G\times\cdots\times G$ ($n$ times)
 has the product structure of a `cubical complex' such that
 each product $\bar c_1\times\cdots\times \bar c_n$ is isometric
 to a Euclidean cube $[0,1]^k$, where $\bar c_i$ is 
 the topological closure of a cell of $G$.
The dimension $k$ is the number of the cells $c_i$ 
 that are edges of $G$.
The \emph{diagonal} of the product $G^n$ is 
$$\De(G^n)=\{(x_1,\dots,x_n)\in G^n\vl x_i=x_j \mbox{ for some }i\neq j\}.$$

\begin{definition}
\label{def:TopConfigurationSpaces}
Let $G$ be a graph, $n$ be a positive integer.
The \emph{ordered topological} configuration space $\OC(G,n)$
 of $n$ distinct robots in $G$ is $G^n-\De(G^n)$.
The \emph{unordered topological} configuration space $\UC(G,n)$
 of $n$ indistinguishable robots in $G$ is the quotient of $\OC(G,n)$
 by the action of the permutation group $S_n$ of $n$ robots.
\end{definition}

The ordered topological space $\OC([0,1],2)$ is the unit square without its diagonal
 $\{(x,y)\in[0,1]^2 \vl x\neq y\}$, which is homotopy equivalent
 to a disjoint union of 2 points.
Topological spaces $X,Y$ are \emph{homotopy} equivalent if 
 there are continuous maps $f:X\to Y$, $g:Y\to X$ such that
 $g\circ f:X\to X$, $f\circ g:Y\to Y$ can be connected 
 with $\id_X:X\to X$, $\id_Y:Y\to Y$, respectively,
 through continuous families of maps.
In particular, $X$ is \emph{contractible} if $X$ is homotopy equivalent to a point.
A space $X$ can be homotopy equivalent to its subspace $Y$ through 
 a \emph{deformation retraction} that is a continuous family of maps 
 $f_t:X\to Y$, $t\in[0,1]$, such that $f_t|_Y=\id_Y$, 
 i.e. all $f_t$ are fixed on $Y$, $f_0=\id_X$ and $f_1(X)=Y$.
\smallskip
 
\begin{figure}[!h]
\includegraphics[scale=1.0]{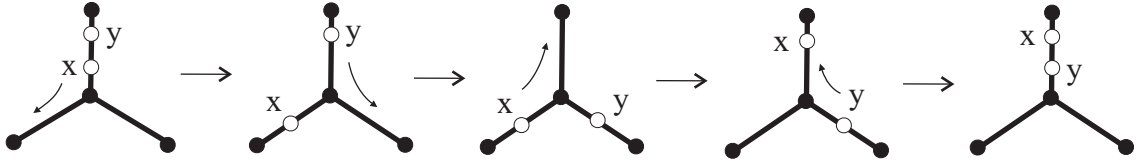}
\caption{Swapping 2 robots $x,y$  without collisions on the triod $T$ }
\label{fig:PermuteRobotsTriod}
\end{figure}

The unordered topological space $\UC([0,1],2)\approx \{(x,y)\in[0,1]^2 \vl x<y\}$
 is contractible to a single point.
More generally, $\OC([0,1],n)$ has $n!$ contractible connected components, while
 $\UC([0,1],n)$ deformation retracts to the standard configuration 
 $x_{i}=(i-1)/(n-1)$, $i=1,\dots,n$, in $[0,1]$.
If a connected graph $G$ has a vertex of degree at least 3 then
 the configuration spaces $\OC(G,n)$, $\UC(G,n)$ are path-connected.
We swap robots $x,y$ near such a vertex as shown in Figure~\ref{fig:PermuteRobotsTriod}.

\begin{definition}
\label{def:GraphBraidGroups}
Given a connected graph $G$ having a vertex of degree at least $3$,
 the \emph{graph braid} groups $\pb(G,n)$ and $\ub(G,n)$ are
 the fundamental groups $\pi_1(\OC(G,n))$ and $\pi_1(\UC(G,n))$,
 respectively, where arbitrary base points are fixed. 
\end{definition}

For the triod $T$ in Figure~\ref{fig:PermuteRobotsTriod},
 both configuration spaces $\OC(T,2)$, $\UC(T,2)$ are
 homotopy equivalent to a circle, see 
 Example~\ref{exa:TopConf2pointT}, i.e. $\ub(T,2)\cong\Z$, 
 $\pb(T,2)\cong\Z$, although $\pb(T,2)$ can be considered
 as an index~2 subgroup $2\Z$ of $\ub(T,2)\cong\Z$.

\begin{definition}
\label{def:DiscConfigurationSpaces}
The \emph{ordered discrete} space $\OD(G,n)$ consists of 
 all the products $\bar c_1\times\cdots\times\bar c_n$ such that 
 each $c_i$ is a cell of $G$ and $\bar c_i\cap\bar c_j=\es$ for $i\neq j$.
The \emph{unordered discrete} space $\UD(G,n)$ is 
 the quotient of $\OD(G,n)$ by the action of $S_n$.
\end{definition}

The \emph{support} $\supp(H)$ of a subset $H\subset G$
 is the minimum union of closed cells containing $H$.
For instance, the support of a vertex or open edge
 coincides with its topological closure in $G$, while 
 the support of a point interior to an open edge $e$ 
 is $\bar e$, i.e. the edge $e$ with its endpoints.
A configuration $(x_1,\dots,x_n)\in G^n$ is \emph{safe} if
 $\supp(x_i)\cap\sup(x_j)=\es$ whenever $i\neq j$.
Then $\OD(G,n)$ consists of all safe configurations:
 $\OD(G,n)=\{(x_1,\dots,x_n)\in G^n\vl 
 \supp(x_i)\cap\supp(x_j)=\es,\; i\neq j\}$.
\smallskip

A path in a graph $G$ is \emph{essential} 
 if it connects distinct essential vertices of $G$.
A cycle in $G$ is \emph{essential} if it contains 
 a vertex of degree more than 2.
Since only connected graphs are considered,
 a non-essential cycle coincides with the whole graph.
Subdivision Theorem~\ref{the:Subdivision} provides sufficient conditions 
 such that the configuration spaces $\OC(G,n),\UC(G,n)$ 
 deformation retract to their discrete analogues $\OD(G,n),\UD(G,n)$, respectively.
Then $\ub(G,n)\cong\pi_1(\UD(G,n))$.

\begin{theorem}
\label{the:Subdivision}
\cite[Theorem~2.1]{Abr1}
Let $G$ be a connected graph, $n\geq 2$.
The discrete spaces $\OD(G,n),\UD(G,n)$ are deformation retracts of 
 the topological configuration spaces $\OC(G,n),UC(g,n)$, respectively, 
 if both conditions (\ref{the:Subdivision}a) and (\ref{the:Subdivision}b) hold:
\smallskip

(\ref{the:Subdivision}a) every essential path in $G$ has at least $n+1$ edges;

(\ref{the:Subdivision}b) every essential cycle in $G$ has at least $n+1$ edges.
\end{theorem}

\begin{figure}[!h]
\includegraphics[scale=1.0]{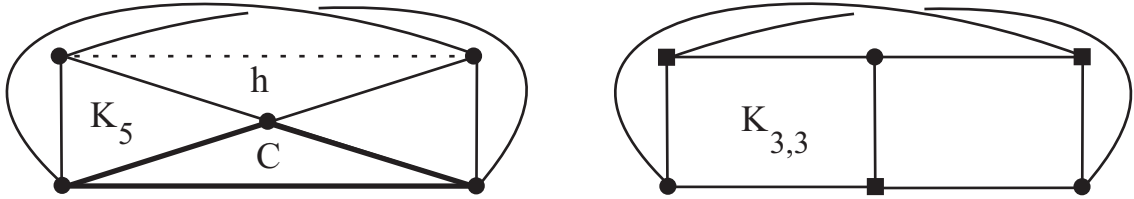}
\caption{Kuratowski graphs $K_5$ and $K_{3,3}$}
\label{fig:KuratowskiGraphs}
\end{figure}

The conditions above imply that $G$ has at least $n$ vertices, so $\OD(G,n)\neq\es$.
A strengthened version of Subdivision Theorem~\ref{the:Subdivision} for $n=2$ 
 only requires that $G$ has no loops and multiple edges \cite[Theorem~2.4]{Abr1}.
Hence the topological configuration spaces of 2 points on 
 the Kuratowski graphs $K_5,K_{3,3}$ deformation retract to 
 their smaller discrete analogues, which are easy to visualise,
 see Figure~\ref{fig:KuratowskiGraphs}.
\smallskip

In $\OD(K_5,2)$, if the 1st robot is moving along an edge $h\in K_5$, then the 2nd robot
 can be only in the triangular cycle $C\subset K_5-h$, which gives in total 
 10 triangular tubes $h\times C$ forming the oriented surface of genus 6.
Similarly, computing the Euler characteristic, we may conclude 
 that  $\OD(K_{3,3},2)$ is the oriented surface of genus 4.
These are the only graphs without loops whose discrete configuration spaces 
 $\OD(G,2)$ are closed manifolds, see \cite[Corollary~5.8]{Abr1}.

\subsection{Main results}
\label{subs:Results}
\noindent
\smallskip

There are two different approaches to computing graph braid groups
 suggested by Abrams \cite[section~3.2]{Abr1} and 
 Farley, Sabalka \cite[Theorem~5.3]{FS}.
In the former approach a graph braid group splits as a graph of simpler groups,
 which gives a nice global structure of the group and proves that, for instance, 
 the graph braid groups are torsion free \cite[Corollary~3.7 on p.~25]{Abr1}.
The latter approach based on the discrete Morse theory by Forman \cite{For}
 writes down presentations of graph braid groups retracting
 a big discrete configuration space to a smaller subcomplex.
\smallskip

We propose another local approach based on classical 
 Seifert -- van Kampen Theorem~\ref{the:SeifertVanKampen}.
Presentations are computed step by step starting from simple graphs and 
 adding edges one by one, which allows us to update growing networks in real-time.
Resulting Algorithm~\ref{alg:GraphBraidGroups} expresses generators 
 of graph braid groups in terms actual motions of robots,
 i.e. as a list of positions at discrete time moments.
We also design motion planning Algorithm~\ref{alg:MotionPlanningUnordered}
 connecting any configurations of $n$ robots.
Its complexity is linear in the number of edges 
 and quadratic in the number of robots.

\begin{algorithm}
\label{alg:GraphBraidGroups}
There is an algorithm writing down a presentation of the graph braid group $\ub(G,n)$ 
 and representing generators by actual paths between configurations of robots, 
 see step-by-step instructions in subsection~\ref{subs:MotionPlanningUnordered}.
\end{algorithm}

\begin{figure}[!h]
\includegraphics[scale=1.0]{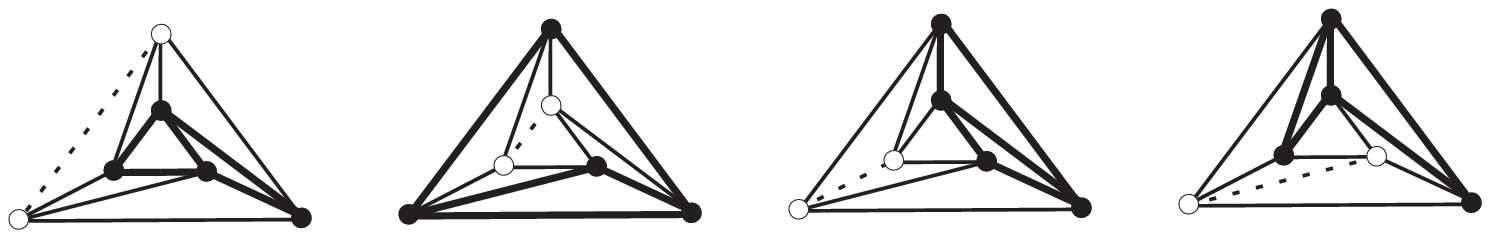}
\caption{A non-light planar graph with 4 choices of a closed edge $\bar h$ }
\label{fig:TriangularGraph}
\end{figure}

According to \cite[Theorem~5.6]{FS2}, the braid groups of planar graphs 
 having only disjoint cycles have presentations where each relator 
 is a commutator, not necessarily a commutator of generators. 
Demonstrating the power of Algorithm~\ref{alg:GraphBraidGroups}, 
 we extend this result to a wider class of light planar graphs.
A planar connected graph $G$ is called \emph{light} if any cycle $C\subset G$ 
 has an open edge $h$ such that all cycles from $G-\bar h$ do not meet $C$.  
Any loop or multiple edge provides an edge $h$ satisfying the above condition.
Figure~\ref{fig:TriangularGraph} shows a non-light planar graph with 
 4 choices of a (dashed) edge $h$ and corresponding (fat) cycles from $G-\bar h$. 
Removing the closure $\bar h$ from $G$ is equivalent to removing 
 the endpoints of $h$ and all open edges attached to them.

\begin{corollary}
\label{cor:2pointGroups}
The braid group $\ub(G,2)$ of any light planar graph $G$ has a presentation 
 where each relator is a commutator of motions along disjoint cycles.
\end{corollary}

A stronger version of Corollary~\ref{cor:2pointGroups} 
 with a geometric description of generators and relators is given in 
 Proposition~\ref{pro:2pointGroupsUnordered}
 in the case of unordered robots.
\medskip

\noindent
{\bf Outline.}
In section~\ref{sect:DiscreteSpaces} we consider 
 basic examples and recall related results.
Section~\ref{sect:FundamentalGroupsUnordered} introduces 
 the engine of Propositions~\ref{pro:AddHangingEdgeUnordered}, 
 \ref{pro:StretchHangingEdgeUnordered},  \ref{pro:CreateCyclesUnordered}  
 updating presentations of graph braid groups by  adding edges one by one.
Section~\ref{sect:ComputingGroupsUnordered} lists step-by-step instructions 
 to compute a presentation of an arbitrary graph braid group.
As an application, we geometrically describe presentations 
 of 2-point braid groups of light planar graphs. 
Further open problems are stated in subsection~\ref{subs:OpenProblems}.
\medskip

\noindent
{\bf Acknowledgements.}
The author thanks Michael Farber for useful discussions and 
 Lucas Sabalka sending an early version of his manuscript \cite{FS2}.
\medskip


\section{Discrete configuration spaces of a graph}
\label{sect:DiscreteSpaces}

In this section we discuss discrete configuration spaces in more details
 and construct them recursively in Lemmas~\ref{lem:RecursiveConstructionUnordered}
 and~\ref{lem:RecursiveConstructionOrdered}. 
Further we assume that $n\geq 2$.


\subsection{Configuration spaces of the triod $T$}
\label{subs:ConfSpacesTriod}
\noindent
\smallskip

\begin{figure}[!h]
\includegraphics[scale=1.0]{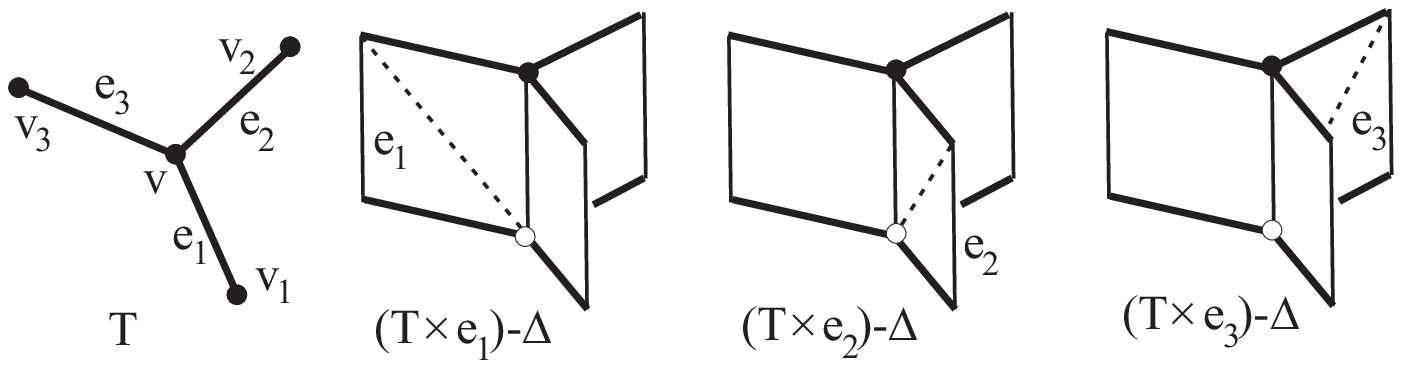}
\caption{The triod $T$ and $(T\times e_1)-\De$, 
 $(T\times e_2)-\De$, $(T\times e_3)-\De$ }
\label{fig:Triod2pointProducts}
\end{figure}

In this subsection we describe configuration spaces of 2 points on the triod $T$
 comprised of 3 hanging edges $e_1,e_2,e_3$ attached to the vertex $v$,
 see Figure~\ref{fig:Triod2pointProducts}.

\begin{figure}[!h]
\includegraphics[scale=1.0]{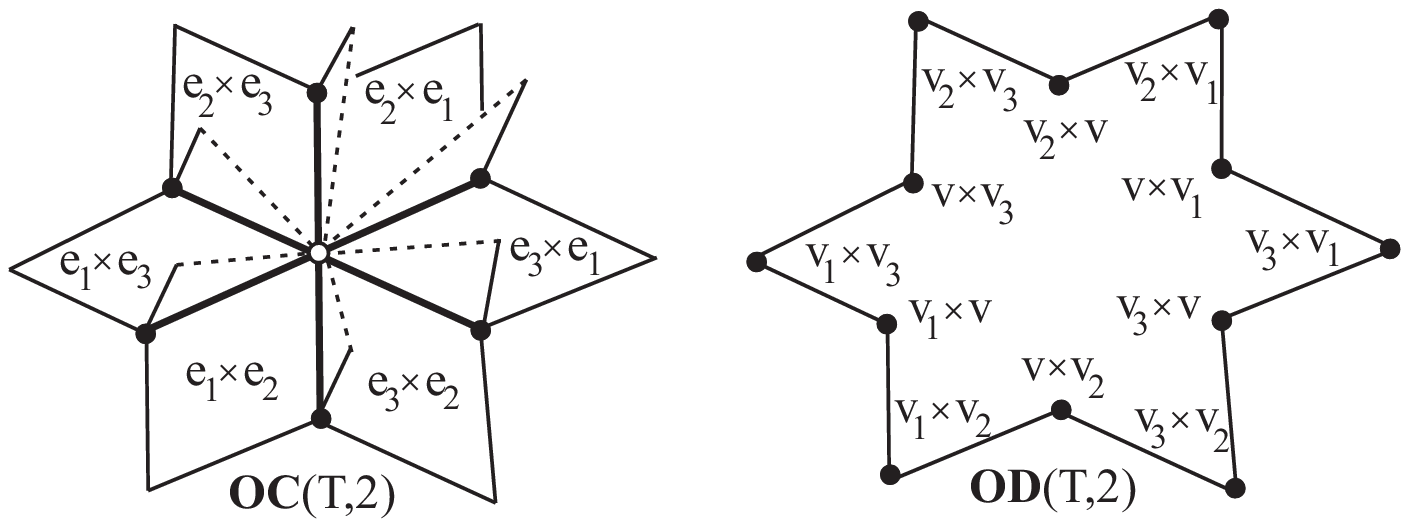}
\caption{The ordered space $\OC(T,2)$ and its discrete analogue $\OD(T,2)$ }
\label{fig:Triod2pointSpaces}
\end{figure}

\begin{example}
\label{exa:TopConf2pointT}
The ordered topological space $\OC(T,2)$ is the union of
 three 3-page books $T\times e_1$, $T\times e_2$, $T\times e_3$
 shown in the right pictures of Figure~\ref{fig:Triod2pointProducts}
 without the diagonal $\De=\{(x,y)\in T^2\vl x\neq y\}$.
Then $\OC(T,2)$ consists of the 6 symmetric rectangles
 $e_i\times e_j$ ($i\neq j$) and 6 triangles from 
 the squares $e_i\times e_i$, $i=1,2,3$, after removing their diagonals, 
 see the left picture of Figure~\ref{fig:Triod2pointSpaces}
 and \cite[Example~6.26]{AF}.
\end{example}

\begin{example}
\label{exa:DiscConf2pointT}
The ordered topological space $\OC(T,2)$ deformation retracts to
 the polygonal circle in the right picture of Figure~\ref{fig:Triod2pointSpaces},
 which is the ordered discrete space $\OD(T,2)$ having 12 vertices
  $v_i\times v_j$ ($i\neq j$) and $v\times v_i$, $v_i\times v$, $i=1,2,3$,
  symmetric under the permutation of factors.
The unordered spaces $\UC(T,2),\UD(T,2)$ are quotients
 of the corresponding ordered spaces by the rotation through $\pi$
  and are homeomorphic to the same spaces $\OC(T,2),\OD(T,2)$, respectively.
Hence the graph braid groups $\ub(T,2)\cong\Z$, $\pb(T,2)\cong\Z$ 
 can be computed using the simpler discrete spaces $\UD(T,2),\OD(T,2)$,
 which is reflected in Subdivision Theorem~\ref{the:Subdivision}.
\end{example}


\subsection{Recursive construction of discrete spaces}
\label{subs:RecursiveConstruction}
\noindent
\smallskip

In this subsection we explain recursive constructions of discrete configuration spaces
 that will be used in section~\ref{sect:FundamentalGroupsUnordered} 
 to compute their fundamental groups.

\begin{figure}[!h]
\includegraphics[scale=1.0]{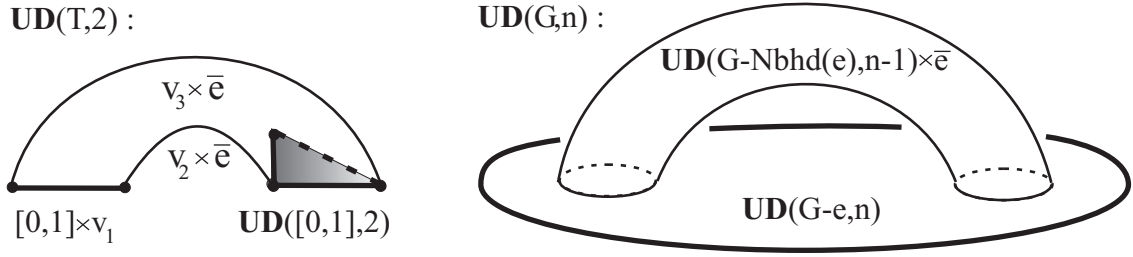}
\caption{Attaching the cylinder in the recursive construction of $\UD(G,n)$ }
\label{fig:AttachCylinder}
\end{figure}

\begin{example}
\label{exa:Recursive2pointT} 
We show how to construct the unordered space 
 $\UD(T,2)$ adding the closed edge $\bar e_1$ to 
 the subgraph $T-(e_1\cup v_1)=\bar e_2\cup\bar e_3\approx [0,1]$.
If both robots $x,y$ are not in the open edge $e_1$, then 
 $(x,y)\in \UD(T-e_1,2)$, where $T-e_1\approx v_1\cup[0,1]$, 
 i.e. either $y=v_1$, $x\in[0,1]$ or $(x,y)\in\UD([0,1],2)$.   
The robot $x$ can not be close to $y$ by Definition~\ref{def:DiscConfigurationSpaces},
 e.g. if $y\in e_1$ then $x\notin e_2\cup e_3$, i.e. $x=v_2$ or $x=v_3$.
Then 
$$\UD(T,2)\approx([0,1]\times v_1)\cup \UD([0,1],2)\cup(\{v_2,v_3\}\times\bar e_1),$$
 where the segments $v_2\times\bar e_1$ and $v_3\times\bar e_1$ are glued 
 at the endpoints $v_2\times v_1,v_3\times v_1$ and $v_2\times v,v_3\times v$,
 respectively.
Up to a homeomorphism, we get 2 arcs attached at theirs endpoints 
 to a solid triangle without one side, see the left picture of 
 Figure~\ref{fig:AttachCylinder}.
\end{example}

The argument of Example~\ref{exa:Recursive2pointT} 
 motivates the following notion.
The \emph{neighbourhood} $\nb(e)$ of an open edge $e\in G$ consists of 
 $\bar e$ and all open edges attached to the endpoints of $e$.
For instance, the complement to the neighbourhood $\nb(e_1)$ 
 in the triod $T$ consists of the hanging vertices $v_2,v_3$, 
 see the left picture of Figure~\ref{fig:Triod2pointProducts}.

\begin{example}
\label{exa:Recursive2pointG} 
Extending the recursive idea of Example~\ref{exa:Recursive2pointT},
 we construct the unordered 2-point space $\UD(G,2)$ of any connected graph $G$.
Fix an open edge $e\subset G$ with vertices $u,v$ and consider the case 
 when one of the robots, say $y$, stays in $e$, then $x\in G-\nb(e)$,
 because $x$ can not be in the same edge $e$ and also in the edges adjacent to $e$.
If both robots $x,y$ are not in $e$ then 
 $(x,y)$ is in the smaller unordered space $\UD(G-e,2)$.
Then $\UD(G,2)$ is a union of smaller subspaces:
 $$\UD(G,2)\approx\UD(G-e,2)\cup ((G-\nb(e))\times \bar e),$$
 where the cylinder $(G-\nb(e))\times \bar e$ is glued to $\UD(G-e,2)$
 along the subgraphs $(G-\nb(e))\times u$ and $(G-\nb(e))\times v$. 
The reduction above extends to a general recursive construction 
 in Lemma~\ref{lem:RecursiveConstructionUnordered}.
Lemmas~\ref{lem:RecursiveConstructionUnordered} and
 \ref{lem:RecursiveConstructionOrdered} are discrete analogues of Ghrist's construction
 of the ordered topological space $\OC(G,n)$ \cite[Lemma~2.1]{Ghr}.
\end{example}

\begin{lemma}
\label{lem:RecursiveConstructionUnordered} 
Let a graph $G$ have an open edge $e$ with vertices $u,v$. 
Then the unordered discrete space $\UD(G,n)$ is homeomorphic to 
 (see Figure~\ref{fig:AttachCylinder})
 $$\UD(G,n)=\UD(G-e,n)\cup (\UD(G-\nb(e),n-1)\times \bar e), \mbox{ where}$$
 the cylinder $\UD(G-\nb(e),n-1)\times \bar e$ 
 is glued to $\UD(G-e,n)$ along 
$$\mbox{the bases }\; \UD(G-\nb(e),n-1)\times u \;
  \mbox{ and }\; \UD(G-\nb(e),n-1)\times v.$$
\end{lemma}
\begin{proof}
In the space $\UD(G,n)$ of all safe configurations $\bx=(x_1,\dots,x_n)$ consider
 the smaller subspace $\UD(G-e,n)$, where $x_i\notin e$ for each $i=1,\dots,n$.
The complement $\UD(G,n)-\UD(G-e,n)$ consists of configurations with(say) $x_n\in e$.
Here the index $n$ is not important since the robots are not ordered.
By Definition~\ref{def:DiscConfigurationSpaces}, 
 the other robots $x_1,\dots,x_{n-1}\notin\nb(e)$, 
 i.e. the complement is 
 $$\UD(G,n)-\UD(G-e,n)\approx\UD(G-\nb(e),n-1)\times e.$$
The bases of the last cylinder are subspaces of the smaller configuration space:
 $$\UD(G-\nb(e),n-1)\times u,\; \UD(G-\nb(e),n-1)\times u\subset \UD(G-e,n).$$
The cylinder $\UD(G-\nb(e),n-1)\times e$ represents motions when the $n$-th robot 
 moves along $e$, while the other robots remain in $\UD(G-\nb(e),n-1)$. 
\end{proof}

Further in sections~\ref{sect:FundamentalGroupsUnordered} 
 and~\ref{sect:ComputingGroupsUnordered}
 the simpler unordered case is considered.
We believe that our approach literally extends to the ordered case 
 using similar Lemma~\ref{lem:RecursiveConstructionOrdered} with 
 $n$ cylinders indexed by $i=1,\dots,n$ since the robots are ordered.

\begin{lemma}
\label{lem:RecursiveConstructionOrdered} 
Let a graph $G$ have an open edge $e$ with vertices $u,v$. 
Then the ordered discrete space $\OD(G,n)$ is homeomorphic to 
 (see Figure~\ref{fig:AttachCylinder})
 $$\OD(G,n)=\OD(G-e,n)\cup_{i=1}^n (\OD_{(i)}(G-\nb(e),n-1)\times\bar e),\mbox{ where}$$
 $\OD_{(i)}(G-\nb(e),n-1)\times\bar e=\{\bx\in\OD(G,n)\vl x_i\in\bar e\}$ 
 is glued to $\OD(G-e,n)$ 
$$\mbox{ along } 
 \OD_{(i)}(G-\nb(e),n-1)\times u=\{\bx\in\OD(G-\nb(e),n)\vl x_i=u\}\mbox{ and}$$
$$\OD_i(G-\nb(e),n-1)\times v=\{\bx\in\OD(G-\nb(e),n)\vl x_i=v\},\; 
  i=1,\dots,n.\eqno{\square}$$
\end{lemma}


\subsection{Homotopy types of configuration spaces}
\label{subs:HomotopyTypes}
\noindent
\smallskip

In this subsection we recall general results on homotopy types of configuration spaces.
Recall that a topological space $X$ is \emph{aspherical} or a $K(\pi,1)$ space
 if it has a contractible universal cover, in particular $\pi_i(X)=0$ for $i>1$.
A covering $p:Y\to X$ is \emph{universal} if the cover $Y$ is simply connected. 
Then the covering $p$ has the \emph{universal} property that, for any covering $q:Z\to X$,
 there is another covering $Y\to Z$ whose composition with $q:Z\to X$ 
 gives the original covering $p:Y\to X$. 

\begin{proposition}
\label{pro:Asphericity} 
\emph{(Asphericity of configuration spaces, 
 Ghrist \cite[Corollary~2.4, Theorem~3.1]{Ghr} for topological spaces 
 and Abrams \cite[section 3.2]{Abr1} for discrete spaces)}
Every component of $\OC(G,n),\UC(G,n),\OD(G,n),\UD(G,n)$ is aspherical.
\qed
\end{proposition}

Ghrist \cite[Corollary~2.4, Theorem~3.1]{Ghr} proves
 the above result for the ordered topological space $\OC(G,n)$,
 which implies the same conclusion for $\UC(G,n)$, because
 the universal cover of a component of $\UC(G,n)$ is
 a universal cover of some component of $\OC(G,n)$
 as mentioned by Abrams \cite[the proof of Corollary~3.6]{Abr1}.
\smallskip

Proposition~\ref{pro:Dimension} implies that the homotopy type of 
 discrete spaces depends on the graph $G$, 
 but not on the number $n$ of robots.
It was proved by Ghrist \cite[Theorems 2.6 and 3.3]{Ghr} 
 for the ordered topological space $\OC(G,n)$,
 which easily extends to the unordered case.
The circle $S^1$ is excluded below, because its unordered space $\UC(S^1,n)$ is
 contractible, while $\OC(S^1,n)$ deformation retracts to a disjoint union of
 $(n-1)!$ configurations indexed by permutations of $n$ robots up to cyclic shifts.

\begin{proposition}
\label{pro:Dimension} 
\emph{(Homotopy type of topological configuration spaces)}
If a connected graph $G$ is not homeomorphic to $S^1$ and has exactly $m$ essential vertices,
 then $\OC(G,n)$ and $\UC(G,n)$ deformation retract to $m$-dimensional complexes.
\qed
\end{proposition}

For instance, the configuration spaces of 2 robots in the triod $T$ 
 having a single essential vertex deformation retract to a 1-dimensonal circle, 
 see Examples~\ref{exa:TopConf2pointT}, \ref{exa:DiscConf2pointT}.


\section{Fundamental groups of unordered discrete spaces}
\label{sect:FundamentalGroupsUnordered}

In this section we compute graph braid groups showing how their presentations 
 change by Seifert -- van Kampen Theorem~\ref{the:SeifertVanKampen} 
 after adding new edges to a graph. 
Let $X,Y$ be open path-connected subsets of $X\cup Y$ 
 such that $X\cap Y\neq\es$ is also path-connected.
If $X,Y$ are not open in $X\cup Y$, they usually can be replaced 
 by their open neighbourhoods that deformation retract to $X,Y$, respectively.
Assume that $X,Y,X\cap Y,X\cup Y$ have a common base point.
If $\ba$ is a finite vector of elements then a group presentation has the form 
 $\ab{\ba\vl\br}$, where the relator $\br$ (a vector of words in the alphabet $\ba$) 
 denotes the vector relation $\br=1$.
We give the practical reformulation of the Seifert -- van Kampen Theorem 
 \cite[Theorem~3.6 on p.~71]{CF}.

\begin{theorem}
\label{the:SeifertVanKampen} 
\emph{(Seifert -- van Kampen Theorem \cite[Theorem~3.6 on p.~71]{CF})}\\
If presentations $\pi_1(X)=\ab{\bb\vl\bla}$, $\pi_1(Y)=\ab{\bg\vl\bmu}$ are given and
 $\pi_1(X\cap Y)$ is generated by (a vector of) words $\ba$, then the group $\pi_1(X\cup Y)$ 
 has the presentation $\pi_1(X\cup Y)=\ab{\bb,\bg\vl \bla,\bmu,\ba_X=\ba_Y}$,
 where $\ba_X,\ba_Y$ are obtained from the words $\ba$ by rewriting 
 them in the alphabets $\bb$, $\bg$, respectively.
\end{theorem}

As an example, consider the 2-dimensional torus $X\cup Y$, where
 $X$ is the complement to a closed disk $D$, while
 $Y$ is a open neighbourhood of $D$, i.e. $X\cap Y$ is an annulus. 
Then $X$ is homotopically equivalent to a wedge of 2 circles, i.e. 
 $\pi_1(X)=\{\al,\be\vl\}$ is free, 
 $\pi_1(Y)=\ab{\vl}$ is trivial and $\pi_1(X\cap Y)=\Z$,
 hence $\pi_1(X\cup Y)=\{\al,\be\vl \al\be\al^{-1}\be^{-1}\}$ as
 $\al\be\al^{-1}\be^{-1}$ represents the boundary of $D$.
\smallskip

We will write down presentations of the fundamental groups 
 $\pi_1(\UD(G,n))\cong\ub(G,n)$ step by step adding edges to the graph 
 and watching the changes in the presentations.
The base of our recursive computation is the contractible space 
 $\UD([0,1],n)$ of $n$ robots in a segment whose fundamental group is trivial.
\smallskip

In Proposition~\ref{pro:AddHangingEdgeUnordered} we glue a hanging edge to
 a vertex of degree at least 2, e.g. to an internal vertex of $[0,1]$, 
 which may create an essential vertex.
In Proposition~\ref{pro:StretchHangingEdgeUnordered} 
 we add a hanging edge to a hanging vertex of degree~1, 
 which does not create an essential vertex.
In Example~\ref{exa:CreateCyclesUnordered} and 
 Proposition~\ref{pro:CreateCyclesUnordered} 
 we attach an edge creating cycles.
Algorithm~\ref{alg:GraphBraidGroups} computing graph braid groups is essentialy based on 
 Propositions~\ref{pro:AddHangingEdgeUnordered}, 
 \ref{pro:StretchHangingEdgeUnordered}, 
 \ref{pro:CreateCyclesUnordered} 
 showing how a presentation is gradually becoming more complicated. 


\subsection{Adding a hanging edge in the unordered case}
\label{subs:AddHangingEdgeUnordered}
\noindent
\smallskip

We start with the degenerate case when a tree $H$ is obtained by 
 adding a hanging edge $e$ to some internal vertex $v$ of $[0,1]$.
Assume that $[0,1]$ is subdivided into at least $n-1$ subedges,
 otherwise the discrete configuration space $\UD(H,n)=\emptyset$
 since $n$ robots occupy at least $n$ distinct vertices.
Choose a hanging (open) edge $e\subset H$ attached to 
 a hanging vertex $u$ and vertex $v$ of degree at least 3.
If the vertex $v$ has degree $\deg v$ then $H-\nb(e)$ consists of 
 $\deg v-1$ disjoint subtrees, some of them could be points.
Hence $\UD(H-\nb(e),n-1)$ splits into $\deg v-1$ subspaces 
 $\UD_j(H-\nb(e),n-1)$, where $j$ may vary from $1$ to $\deg v-1$.
Fix base points: 
$$a\in\UD(H-(e\cup u),n),\quad
  c_j\in\UD_{j}(H-\nb(e),n-1).$$

\begin{figure}[!h]
\includegraphics[scale=1.0]{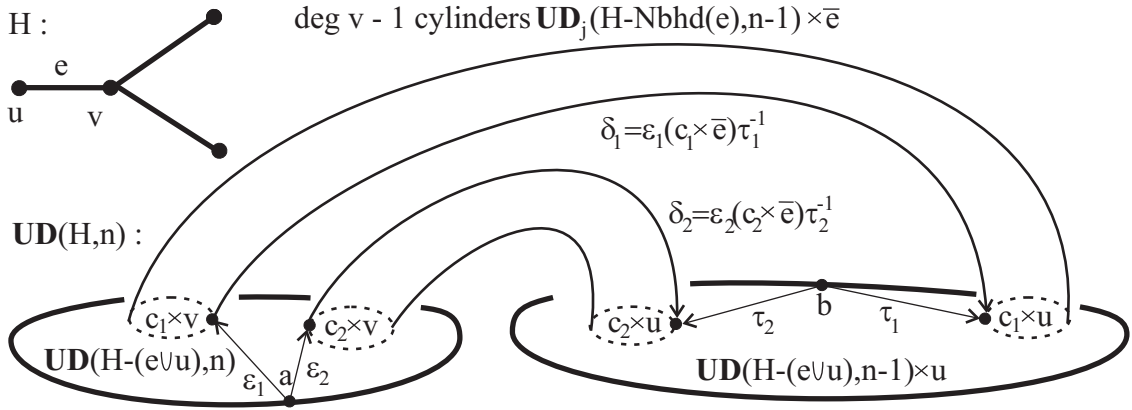}
\caption{Adding a hanging edge $e$ to a non-hanging vertex $v$ }
\label{fig:AddHangingEdge}
\end{figure}

We also fix a base point $b\in\UD(H-(e\cup u),n-1)\times u$,
 which can be chosen as $c_1\times u$ for simplicity.
In $\UD(H-\nb(e),n-1)$ find a path $\e_j$ from $a$ to $c_j\times v$,
 a path $\tau_j$ from $b$ to $c_j\times u$, $j=1,\dots,\deg v-1$,
 see Figure~\ref{fig:AddHangingEdge} and 
 motion planning  algorithm~\ref{alg:MotionPlanningUnordered}
 in subsection~\ref{subs:MotionPlanningUnordered}.
The base configurations $a_j,b_j$ are connected by the motion $(c_j\times\bar e)$
 when $n-1$ robots stay fixed at $c_j\in\UD(H-\nb(e),n-1)$ and 
 1 robot moves along $\bar e$, see Figure~\ref{fig:AddHangingEdge}.
Adding $\e_j,\tau_j^{-1}$ at the start and end of the motion  $(c_j\times\bar e)$, 
 respectively, we get the $\deg v-1$ paths $\de_j$ going from $a$ to $b$ 
 in $\UD(H,n)$, $j=1,\dots,\deg v-1$.
For a loop $\be\subset\UD(H-\nb(e),n-1)$ representing a motion of $n-1$ robots,
 the loop $(\be\{x_n=u\})\subset\UD(H-\nb(e),n-1)\times u$ denotes the motion
 when $n-1$ robots follow $\be$ and one robot remains fixed at $u$.

\begin{proposition}
\label{pro:AddHangingEdgeUnordered} 
\emph{(Adding a hanging edge $e$ to a non-hanging vertex $v$)}\\
In the notations above
 and for presentations $\pi_1(\UD(H-(e\cup u),n))=\ab{ \ba \vl \br }$ and
$$\pi_1(\UD(H-(e\cup u),n-1))=\ab{ \bb \vl \bla },\;
  \pi_1(\UD_j(H-\nb(e),n-1))=\ab{ \bg_j \vl \bmu_j },$$
 the group $\pi_1(\UD(H,n))$ is generated by 
 $\quad\ba$, $\quad\de_1(\bb\{x_n=u\})\de_1^{-1}$, $\quad\de_1\de_j^{-1} (j>1)$, 
$$\mbox{ subject to }\quad
 \br=1,\quad \de_1(\bla\{x_n=u\})\de_1^{-1}=1,\quad
 (\bg_j\{x_n=v\})=\de_j(\bg_j\{x_n=u\})\de_j^{-1}.$$
\end{proposition}
\begin{proof}
By the recursive construction from Lemma~\ref{lem:RecursiveConstructionUnordered} 
 one has
 $$\UD(H,n)\approx\UD(H-e,n)\cup (\UD(H-\nb(e),n-1)\times\bar e).$$
Since $H-e$ splits into the vertex $u$ and the remaining subgraph $H-(e\cup u)$,
 then the space $\UD(H-e,n)$ consists of the 2 connected components
 $\UD(H-(e\cup u),n)$, where all robots are in $H-(e\cup u)$, and 
 $\UD(H-(e\cup u),n-1)\times u$, where one robot is at $u$.
The non-connected cylinder $\UD(H-\nb(e),n-1)\times\bar e$ splits
 into $\deg v-1$ cylinders $\UD_j(H-\nb(e),n-1)\times\bar e$
 connecting $\UD(H-(e\cup u),n)$ and $\UD(H-(e\cup u),n-1)\times u$
 since the complement $H-\nb(e)$ is obtained from $H$ by removing $u,v$ 
 and all open edges attached to the vertex $v$ of degree $\deg v$.
\smallskip

Add the cylinders $\UD_j(H-\nb(e),n-1)\times\bar e$
 to the subspace $\UD(H-(e\cup u),n)$, which
 does not affect the group $\pi_1(\UD(H-(e\cup u),n))$, because
 the cylinders deformation retract to their bases $\UD_j(H-\nb(e),n-1)\times v$.
To apply Seifert -- van Kampen Theorem~\ref{the:SeifertVanKampen} correctly, 
 add all the paths $\de_j$ to the resulting union, which gives
 the $\deg v-2$ new generators $\de_1\de_j^{-1}$, $j>1$.
\smallskip

Consider the space $\UD(H-(e\cup u),n-1)\times u$ as a subspace of $\UD(H,n)$.
Formally a loop $\be\in\pi_1(\UD(H-(e\cup u),n-1))$ becomes
 the loop $(\be\{x_n=u\})$ from $\pi_1(\UD(H-(e\cup u),n-1)\times u)$,
 where one robot remains fixed at $u$.
The same argument applies to the relator $\bla$.
No other relations appear as the intersection of $\cup_j\de_j$ and
 $\UD(H-(e\cup u),n)\cup_j(\UD_j(H-\nb(e),n-1)\times\bar e)$ contracts to $a$.
\smallskip

Now take the union with the remaining subspace $\UD(H-(e\cup u),n-1)\times u$, which 
 adds the generators and relations of $\pi_1(\UD(H-(e\cup u),n-1))=\ab{ \bb \vl \bla }$.
The resulting intersection deformation retracts to the wedge of 
 the $\deg v-1$ bases $\UD_j(H-\nb(e),n-1)\times u$, 
 so each generator $\bg_j$ gives a relation between
 the words representing the loops $(\bg_j\{x_n=v\})$ in 
 the spaces $\UD(H-(e\cup u),n)$ and $\UD(H-(e\cup u),n-1)\times u$.
In the latter space the loop can be conjugated by $\de_j$,
 which replaces $b$ by the base point $a\in \UD(H,n)$, we may set $j=1$.
\smallskip

Notice that the loops $\de_j(\bg_j\{x_n=u\})\de_j^{-1}$ live in 
 $\UD(H-(e\cup u),n)$ with the base point $a$ and can be expressed 
 in terms of the generators $\de_j(\bb\{x_n=u\})\de_j^{-1}$.
So the last equality in the presentation is 
 a valid relation between new generators.
\end{proof}


\subsection{Stretching a hanging edge in the unordered case}
\label{subs:StretchHangingEdgeUnordered}
\noindent
\smallskip

\begin{figure}[!h]
\includegraphics[scale=1.0]{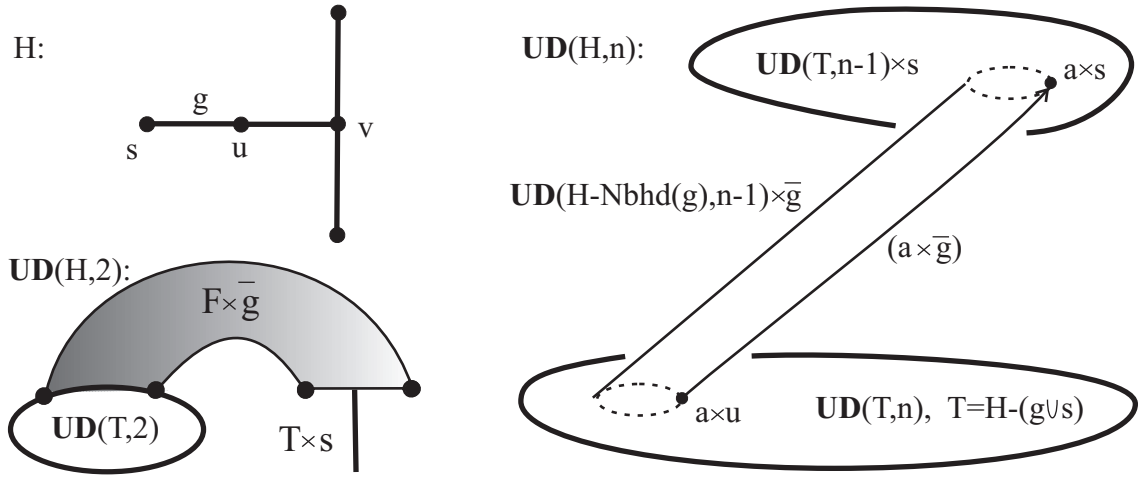}
\caption{Stretching a hanging edge in a tree $H$}
\label{fig:StretchHangingEdge}
\end{figure}

In this subsection we show how the presentation of a braid group 
 changes after stretching a hanging edge of a tree.
First we consider the degenerate case of stretching a hanging edge $e$ of 
 the triod $T$ in the top left picture of  Figure~\ref{fig:StretchHangingEdge}. 

\begin{example}
\label{exa:StretchHangingEdgeUnordered}
Let $H$ be the tree obtained by adding a hanging edge $g$ 
 to the hanging vertex $u$ of the triod $T$
 in the top left picture of Figure~\ref{fig:StretchHangingEdge}, i.e.
 $T=H-(g\cup s)$, where $s$ is the only hanging vertex of $g$ in the tree $H$.
The complement $F=H-\nb(g)$ consists of 2 hanging edges 
 distinct from $e$ and meeting at the centre $v$ of the triod $T$.
We compute the braid group $\ub(H,2)$ using 
 $\ub(T,2)\cong\Z$ from Example~\ref{exa:DiscConf2pointT}.  
By Lemma~\ref{lem:RecursiveConstructionUnordered} 
 the unordered space $\UD(H,2)$ has the form
 $$\UD(H,2)\approx\UD(H-g,2)\cup(F\times\bar g)=
 \UD(T,2)\cup(T\times s)\cup(F\times\bar g),$$
 where the 2 components of $\UD(H-g,2)$ are connected by the band $F\times\bar g$.
First we apply Seifert -- van Kampen Theorem~\ref{the:SeifertVanKampen} 
 to the union $\UD(T,2)\cup(F\times\bar g)$, which keeps the fundamental group 
 unchanged, i.e. isomorphic to $\ub(T,2)\cong\Z$, because 
 the union deformation retracts to $\UD(T,2)$.
Then we apply the same trick taking the union with 
 $T\times s$, which leads to $\ub(H,2)\cong\Z$ for the same reasons.
\end{example}

Proposition~\ref{pro:StretchHangingEdgeUnordered} below extends 
 Example~\ref{exa:StretchHangingEdgeUnordered} to a general tree $H$.
Choose an (open) edge $g\subset H$ with 
 a hanging vertex $s$ and vertex $u$ of degree~2. 
Fix a base point:
$$a\in\UD(H-\nb(g),n-1)\subset\UD(H-(g\cup s),n-1).$$
Let $(a\times\bar g)$ be the motion from $a\times u$ to $a\times s$ in $\UD(H,n)$,
 when $n-1$ robots stay fixed at $a$, while 1 robot moves along $\bar g$,
 see the right picture of Figure~\ref{fig:StretchHangingEdge}.
Then, for a loop $\ga\in\pi_1(\UD(H-\nb(g),n-1))$, both loops 
 $(\ga\{x_n=u\})$ and $(a\times\bar g)^{-1}(\ga\{x_n=s\})(a\times\bar g)$ 
 pass through the base point $a\times u\in \UD(H,n)$.
\smallskip

\begin{proposition}
\label{pro:StretchHangingEdgeUnordered} 
\emph{(Stretching a hanging edge)}\\
In the notations above and for presentations
$\pi_1(\UD(H-(g\cup s),n))=\ab{ \ba \vl \br }$ and
$$\pi_1(\UD(H-(g\cup s),n-1))=\ab{ \bb \vl \bla },\;
    \pi_1(\UD(H-\nb(g),n-1))=\ab{ \bg \vl \bmu },$$
$\pi_1(\UD(H,n)) \mbox{ is generated by } \ba,\; 
 (a\times\bar g)(\bb\{x_n=s\})(a\times\bar g)^{-1}
 \mbox{ subject to }  \br=1,$
$$ (a\times\bar g)(\bla\{x_n=s\})(a\times\bar g)^{-1}=1,\; 
(\bg\{x_n=u\})=(a\times\bar g)(\bg\{x_n=s\})(a\times\bar g)^{-1}.$$
\end{proposition}
\begin{proof}
By the recursive construction from Lemma~\ref{lem:RecursiveConstructionUnordered} 
 one has
 $$\UD(H,n)\approx\UD(H-g,n)\cup (\UD(H-\nb(g),n-1)\times\bar g),$$
 where the cylinder $\UD(H-\nb(g),n-1)\times\bar e$ 
 is glued to $\UD(H-g,n)$ along the bases $\UD(H-\nb(g),n-1)\times s$ 
 and $\UD(H-\nb(g),n-1)\times u$.
Since $g$ is hanging then $H-\nb(g)$ has 2 components:
 the hanging vertex $s$ and remaining tree $T=H-(g\cup s)$,
 hence $\UD(H-g,n)\approx\UD(T,n)\cup(\UD(T,n-1)\times s)$.
\smallskip

Since the edge $e$ is hanging in $H-(g\cup s)$ before stretching then the complement 
 $H-\nb(g)$ and cylinder $\UD(H-\nb(g),n-1)\times\bar g$ are connected.
Adding the cylinder to $\UD(T,n)$ does not change the presentation of 
 the fundamental group, because the cylinder deformation retracts 
 to its base in $\UD(T,n)$.
Then add $\UD(T,n-1)\times s$ meeting
 the previous union along $\UD(H-\nb(g),n-1)\times s$.
\smallskip

By Seifert -- van Kampen Theorem~\ref{the:SeifertVanKampen} 
 to get a presentation of $\pi_1(\UD(H,n))$ with the base point $a\times u$, 
 we add the generators $(a\times\bar g)(\bb\{x_n=s\})(a\times\bar g)^{-1}$  
 and relations $(a\times\bar g)(\bla\{x_n=s\})(a\times\bar g)^{-1}$ 
 coming from  the group $\pi_1(\UD(T,n-1))$.
Add the new relations $(\bg\{x_n=u\})=(a\times\bar g)(\bg\{x_n=s\})(a\times\bar g)^{-1}$ 
 saying that the generators of the group $\pi_1(\UD(H-\nb(g),n-1))$ 
 after adding the stationary $n$-th robot  become homotopic 
 through the subspace $\UD(H-\nb(g),n-1)\times\bar g$.
\end{proof}


\subsection{Creating cycles in the unordered case}
\label{subs:CreateCyclesUnordered}
\smallskip
\noindent

In this subsection we extend our computations to graphs containing cycles.
First we show how the braid group changes if
 an edge is added at 2 vertices of a triod. 

\begin{figure}[!h]
\includegraphics[scale=1.0]{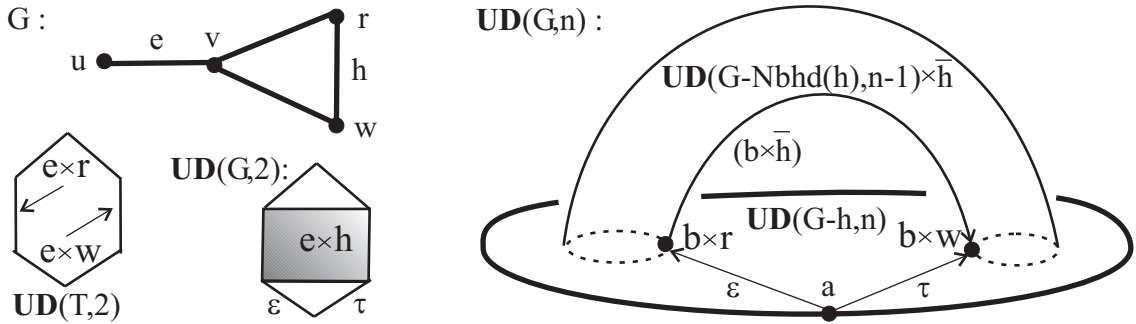}
\caption{Adding an edge $h$ creating cycles }
\label{fig:CreateCycles}
\end{figure}

\begin{example}
\label{exa:CreateCyclesUnordered} 
Let $G$ be the graph obtained from the triod $T$ in the top left picture of 
 Figure~\ref{fig:CreateCycles} by adding the edge $h$ at the vertices $r,w$.
By Lemma~\ref{lem:RecursiveConstructionUnordered} one has
 $$\UD(G,2)\approx\UD(G-h,2)\cup((G-\nb(h))\times\bar e)
 \approx\UD(T,2)\cup(\bar e\times\bar h).$$
Geometrically the band $\bar e\times\bar h$ is glued to the hexagon $\UD(T,2)$
 as shown in the bottom left picture of Figure~\ref{fig:CreateCycles}.
To compute the graph braid group $\ub(G,2)$
 we first add to the band $\bar e\times\bar h$ the motions $\e,\tau\subset\UD(T,2)$ 
 connecting the base configuration $u\times v$ to $u\times r$, $u\times w$, respectively. 
This adds a generator to the trivial fundamental group 
 of the contractible band $\bar e\times\bar h$.
Second we add the union $(\bar e\times\bar h)\cup (\e\cup\tau)$ 
 to $\UD(T,2)$, which gives $\UD(G,2)$.
The intersection of the spaces attached above has the form 
 $(\bar e\times r)\cup(u\times\bar h)\cup(\bar e\times w)$  and is contractible, i.e. 
 $\ub(G,2)$ is the free product of $\ub(T,2)=\Z$ and 
 $\pi_1((\bar e\times\bar h)\cup \e\cup\tau)=\Z$.
\end{example}

Proposition~\ref{pro:CreateCyclesUnordered} extends 
 Example~\ref{exa:CreateCyclesUnordered}
 to a general graph excluding the case $G\approx S^1$.
Choose an (open) edge $h\subset G$ with vertices $r,w$ such that $G-h$ is connected.
Let $G-\nb(h)$ consist of $k$ connected components.
Then $\UD(G-\nb(h),n-1)$ splits into 
 $k$ subspaces $\UD_j(G-\nb(h),n-1)$, where $j=1,\dots,k$.
Fix base points $a\in\UD(G-h,n)$ and $b_j\in\UD_j(G-\nb(h),n-1)$.
Denote by $(b_j\times h)\subset\UD(G,n)$ the motion such that
 one robot goes along the path $(b_j\times h)$ from $b_j\times r$ to $b_j\times w$, 
 while the other robots remain fixed at $b_j\in\UD_j(G-\nb(h),n-1)$, 
 see the right picture of Figure~\ref{fig:CreateCycles} 
 in the case $k=1$ when we may skip the index $j$.
Take paths $\e_j,\tau_j$ going from $a$ to $b_j\times r,b_j\times w$, 
 respectively, in $\UD(G-h,n)$, see Algorithm~\ref{alg:MotionPlanningUnordered}.
Then $\e_j(b_j\times h)\tau_j^{-1}$ is a loop with the base point $a$ 
 in the space $\UD(G,n)$.

\begin{proposition}
\label{pro:CreateCyclesUnordered} 
\emph{(Adding an edge $h$ creating cycles)}
Given presentations
$$\pi_1(\UD(G-h,n))=\ab{ \ba \vl \br },\;
   \pi_1(\UD_j(G-\nb(h),n-1))=\ab{ \bb_j \vl \bla_j }, j=1,\dots,k,$$
 the group $\pi_1(\UD(G,n))$ is generated by $\ba$, $\e_j(b_j\times h)\tau_j^{-1}$
 subject to $\br=1$ and 
$$ \e_j(\bb_j\{x_n=r\})\e_j^{-1}=(\e_j(b_j\times h)\tau_j^{-1})\cdot
 (\tau_j(\bb_j\{x_n=w\})\tau_j^{-1})\cdot(\e_j(b_j\times h)\tau_j^{-1})^{-1}.$$
\end{proposition}
\begin{proof}
The $k$ subspaces $\UD_j(G-\nb(h),n-1)$ can be disconnected, but they are in 
 a 1-1 correspondence with the connected components of $G-\nb(h)$.
Each of the cylinders $\UD_j(G-\nb(h),n-1)\times\bar h$ meets 
 the subspace $\UD(G-h,n)$ at the bases $\UD_j(G-\nb(h),n-1)\times r$ 
 and $\UD_j(G-\nb(h),n-1)\times w$.
\smallskip

First we add to each cylinder $\UD_j(G-\nb(h),n-1)\times\bar h$ the union of 
 the paths $\e_j\cup \tau_j$ connecting the bases to $a\in \UD(G-h,n)$, 
 see Figure~\ref{fig:CreateCycles}.
The fundamental group of $(\UD_j(G-\nb(h),n-1)\times\bar h)\cup (\e_j\cup \tau_j)$ is isomorphic to
 the free product of $\ub(G-\nb(h),n-1)$ and $\Z$ generated by the loop $\e_j(b_j\times h)\tau_j^{-1}$.
Second we add to $\UD(G-h,n)$ each union $(\UD_j(G-\nb(h),n-1)\times\bar h)\cup (\e_j\cup \tau_j)$.
The intersection of the spaces attached above has the form
$$(\UD_j(G-\nb(h),n-1)\times r)\cup (\e_j\cup \tau_j)\cup 
 (\UD_j(G-\nb(h),n-1)\times w)$$
 and is homotopically a wedge of 2 copies of the base $\UD_j(G-\nb(h),n-1)$.
By Siefert -- van Kampen Theorem~\ref{the:SeifertVanKampen} we express 
 the loops $\e_j(\bb_j\{x_n=r\})\e_j^{-1}$ and $\tau_j(\bb_j\{x_n=w\})\tau_j^{-1}$ 
 generating the fundamental group of the intersection in terms of the loops from  
 $\UD(G-h,n)\mbox{ and } (\UD_j(G-\nb(h),n-1)\times\bar h)\cup (\e_j\cup \tau_j).$
In the latter space these loops are conjugated by 
 $\e_j(b_j\times h)\tau_j^{-1}$ as required,
 i.e. homotopic through the cylinder $\UD_j(G-\nb(h),n-1)\times\bar h$.
\end{proof}

If the vector of generators $\br$ is empty, i.e. the groups 
 $\pi_1(\UD_j(G-\nb(h),n-1))$ are trivial, then no new relations 
 are added in Proposition~\ref{pro:CreateCyclesUnordered}.


\section{Computing graph braid groups}
\label{sect:ComputingGroupsUnordered}

At the end of subsection~\ref{subs:MotionPlanningUnordered}
 we give step-by-step instructions of Algorithm~\ref{alg:GraphBraidGroups}
 computing presentations of graph braid groups.
The computing algorithm is based on the technical propositions from 
 section~\ref{sect:FundamentalGroupsUnordered} and auxiliary algorithms 
 from subsection~\ref{subs:MotionPlanningUnordered} below.
As a theoretical application, in Proposition~\ref{pro:2pointGroupsUnordered} we extend  
 the result about 2-point braid groups of graphs with only disjoint cycles \cite[Theorem~5.6]{FS2}
 to a wider class of graphs including all light planar graphs.


\subsection{A motion planning algorithm}
\label{subs:MotionPlanningUnordered}
\noindent
\smallskip

Proposition~\ref{pro:AddHangingEdgeUnordered} requires
 a collision free motion connecting two configurations of $n$ robots.
Take a connected graph $G$ and number its vertices.
We will work with discrete configuration spaces assuming that at every discrete moment 
 all robots are at vertices of a graph $G$ and in one step any robot 
 can move to an adjacent vertex if it is not occupied.
The output contains positions of all robots at every moment.
\smallskip

To describe planning Algorithm~\ref{alg:MotionPlanningUnordered} 
 we introduce auxiliary definitions and searching 
 Algorithms~\ref{alg:ExtremeRobots}, \ref{alg:NeighbourRobot}. 
The $i$-th robot is called \emph{extreme} in a given configuration $(x_1,\dots,x_n)\in\UD(G,n)$ 
 if the remaining robots are in one connected component of $G-x_i$.
One configuration may have several extreme robots, e.g. on a segment 
 there are always 2 extreme robots, while on a circle every robot is extreme.

\begin{algorithm}
\label{alg:ExtremeRobots}
If a graph $G$ has $l$ edges then there is an algorithm of complexity $O(nl)$
 finding all extreme robots in a configuration $(x_1,\dots,x_n)\in\UD(G,n)$.
\end{algorithm}
\begin{proof}
For each robot $x_i$ we visit all vertices of $G-x_i$ remembering the robots we have seen.
If not all robots were seen then the robot $x_i$ is not extreme and we check a robot 
 from a smaller connected component of $G-x_i$, which has fewer edges than $G$.
Hence we will inevitably find an extreme robot, which requires 
 in total not more than $l$ steps for each $i=1,\dots,n$.
\end{proof}

A robot $x_j$ is a \emph{neighbour} of a robot $x_i$ if
 a shortest path from $x_j$ to $x_i$ has the minimal number of edges 
 among all shortest paths from $x_j$ to robots $x_k$ for $k\neq i$.
For $n$ robots on a segment each of the 2 extreme robots has
 a unique neighbour, while on a circle each robot has 2 neighbours.
A shortest path to a neighbour does not contain other robots,
 i.e. the corresponding motion is collision free. 

\begin{algorithm}
\label{alg:NeighbourRobot}
If a connected graph $G$ has $l$ edges then there is an algorithm of 
 complexity $O(l)$ finding a shortest path from a robot $x_i$ to its neighbour $x_j$. 
\end{algorithm}
\begin{proof}
We travel on $G$ in a `spiral way' starting from $x_i$,
 i.e. first we visit all vertices adjacent to $x_i$ and check if there is 
 another robot $x_j$ at one of them, which can be a neighbour of $x_i$.
If not then repeat the same procedure recursively for all these adjacent vertices.
In total we pass through not more than $l$ edges of $G$.
\end{proof}

\begin{algorithm}
\label{alg:MotionPlanningUnordered}
If a connected graph $G$ has $l$ edges, there is an algorithm of  complexity $O(n^2l)$ 
 finding a motion between configurations of $n$ robots in $\UD(G,n)$.
\end{algorithm}
\begin{proof}
For simplicity we assume that all robots are at vertices of degree~2,
 otherwise we may subdivide edges of the graph $G$ and 
 move a robot to an adjacent vertex of degree 2.
This increases the number $l$ of edges by not more than $n\leq l$.
\smallskip

\noindent
\emph{Step 1.} 
Using Algorithm~\ref{alg:ExtremeRobots} of complexity $O(nl)$,
 find an extreme robot in the collection of $2n$ given positions (initial and final together).
\smallskip

\noindent
\emph{Step 2.} 
Assume that the found extreme robot, say $y_n$, is from 
 the final configuration, otherwise swap the roles of initial and final positions.
Using Algorithm~\ref{alg:ExtremeRobots} of complexity $O(l)$, find a shortest path 
 from $y_n$ to its neighbour, say $x_n$, from the initial configuration.
Then safely move $x_n$ towards $y_n$ along the shortest path avoiding 
 collisions and keeping fixed all other robots from the initial configuration.
\smallskip

\noindent
\emph{Step 3.}
Remove from the graph $G$ the robot $y_n$ at a vertex of degree~2
 and all open edges attached to $y_n$ reducing the problem 
 to a smaller graph with $n-1$ robots.
The new graph remains connected since the robot $y_n$ was extreme. 
Return to \emph{Step 1} applying the recursion $n-1$ times,
 which gives $O(n^2l)$ operations in total.
\end{proof}

In Algorithm~\ref{alg:MotionPlanningUnordered} the quadratic complexity 
 in the number of robots seems to be asymptotically optimal, because
 avoiding collisions between $n$ robots should involve 
 some analysis of their pairwise positions. 
\smallskip

\noindent
{\bf Step-by-step instructions of Algorithm~\ref{alg:GraphBraidGroups}.}

\noindent
Start from $n$ robots on a segment subdivided into $n-1$ subsegments, when 
 the configuration space $\UD([0,1],n)$ is a single point and $\ub([0,1],n)$ is trivial.
Construct the graph $G$ adding edges one by one and
 updating presentations of resulting graph braid groups by
 Propositions~\ref{pro:AddHangingEdgeUnordered}, \ref{pro:StretchHangingEdgeUnordered}
 and \ref{pro:CreateCyclesUnordered}.
When we need a motion connecting 2 configurations, we apply 
 motion planning Algorithm~\ref{alg:MotionPlanningUnordered}.
Every generator is represented as a list of vertices 
 where robots are located at every discrete moment. 
 

\subsection{2-point braid groups of graphs in the unordered case}
\label{subs:2pointBraidGroupsPlanar}
\noindent
\smallskip

The first part of Lemma~\ref{lem:2pointTreeGroupsUnordered} without computing 
 the rank was obtained by the global approach of Abrams \cite[Corollary~]{Abr1}.
The second part was claimed by Farber \cite[Theorems~9, 10]{Far1}.
Both parts follow from our local step-by-step computations.

\begin{lemma}
\label{lem:2pointTreeGroupsUnordered} 
For any tree $H$, the braid group $\ub(H,2)$ is free
 and has the rank\\ $\sum(\deg v -1)(\deg v-2)/2$, where
 the sum is over all vertices of degree at least 3.
\end{lemma}
\begin{proof}
Induction on the number of edges of $H$.
The base $H\approx[0,1]$ is trivial.
In the inductive step notice that trees are contractible,
 hence their fundamental groups are trivial and for $n=2$
 the vectors $\bb,\bg,\bla,\bmu$ (with indices $j$) are empty 
 in Propositions~\ref{pro:AddHangingEdgeUnordered} 
 and~\ref{pro:StretchHangingEdgeUnordered}.
The vectors $\br$ are also empty, because they can 
 only come from 2-point braid groups of smaller trees.
So the braid group $\ub(H,2)$ is free.
The only generators of $\ub(H,2)$ are $\de_1\de_j^{-1}$, $j=2,\dots,\deg v-1$, 
 coming from Proposition~\ref{pro:AddHangingEdgeUnordered}, which gives 
 $1+2+\dots+(\deg v-2)=(\deg v -1)(\deg v -2)/2$ generators 
 in total after attaching all edges to each vertex $v$ of degree $\deg v$.
\end{proof}

The Kuratowski graphs $K_5,K_{3,3}$ in Figure~\ref{fig:KuratowskiGraphs} do not satisfy 
 Lemma~\ref{lem:ChooseEdge}, because  the complement to the neighbourhood 
 of any edge $h\in K_5$ ($h\in K_{3,3}$, respectively)  is the triangular 
 (rectangular, respectively) cycle  intersecting any cycle $C\supset h$.

\begin{lemma}
\label{lem:ChooseEdge}
Any light planar graph can be constructed from a tree by adding edges as follows:
 an open edge $h$ added to the new graph $G$ creates a cycle $C$ not meeting 
 any cycle from $G-\nb(h)$ having all its cycles in one connected component.
\end{lemma}
\begin{proof}
Recall that a planar connected graph $G$ is light if any cycle $C\subset G$ has an edge $h$
 such that all cycles from $G-\bar h$ (or, equivalently, $G-\nb(h)$) do not meet $C$.
For a given light planar graph $G$, take any cycle $C$ and corresponding edge $h$.
The smaller graph $G-h$ is light planar, because 
 it has fewer cycles satisfying the same condition.
We may also assume that all cycles of the subgraph $G-\nb(h)$ are 
 in one connected component, otherwise it splits as in 
 the left picture of Figure~\ref{fig:ChooseEdge}.

\begin{figure}[!h]
\includegraphics[scale=1.0]{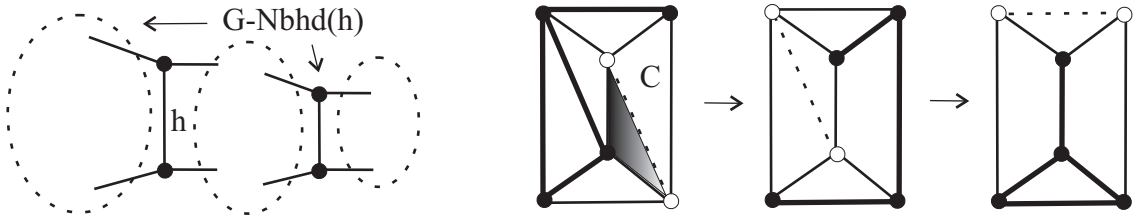}
\caption{Choosing an edge $h$ and a cycle $C\supset h$ in Lemma~\ref{lem:ChooseEdge} }
\label{fig:ChooseEdge}
\end{figure}

Indeed, the open edge $h$ can not split $G$ since  $h$ belongs to the cycle $C\subset G$.
Then we may choose another cycle from a component of $G-\nb(h)$ 
 with a smaller number of edges etc.
Remove edges one by one until the light planar graph becomes a tree.
The original graph can be reconstructed by reversing the procedure above.
\end{proof}

The construction from Lemma~\ref{lem:ChooseEdge} 
 is also applicable to some non-light planar graphs.
The right picture of Figure~\ref{fig:ChooseEdge} shows 3 stages of such a construction, where 
 the closed edge $\bar h$ is dashed and the corresponding subgraph $G-\nb(h)$ has fat edges.
The biggest graph fails to be light planar because of the cycle bounding the grey triangle. 
For the same graph and dashed edge $h$, one can choose another cycle $C$
 that does not meet the only (triangular) cycle from $G-\nb(h)$.
Lemma~\ref{lem:ChooseEdge} implies that Corollary~\ref{cor:2pointGroups} 
 for unordered robots is a particular case of more technical 
 Proposition~\ref{pro:2pointGroupsUnordered}, which holds
 for all graphs constructed as described above.

\begin{proposition}
\label{pro:2pointGroupsUnordered}
For any graph $G$ constructed from a tree as in 
 Lemma~\ref{lem:ChooseEdge}, let $m$ be the first Betti number of $G$.
The braid group $\ub(G,2)$ has a presentation with
 $m+\sum(\deg v-1)(\deg-2)/2$ generators subject to commutator relations, 
 where the sum is over all vertices $v\in G$ of degree at least 3.
A geometric description follows.
\smallskip

\noindent
$\bullet$
At each vertex $v\in G$ fix an edge $e_0$.
For any unordered pair of other edges $e_i,e_j$ at the same vertex $v$, $j=1,\dots,\deg v-1$,
 one generator of $\ub(G,2)$ swaps 2 robots in the triod $e_0\cup e_i\cup e_j$ 
 using the collision free motion shown in Figure~\ref{fig:PermuteRobotsTriod}.
\smallskip

\noindent
$\bullet$
Denote by $h_1,\dots,h_m$ disjoint open edges of $G$ such that $G-(\cup_{j=1}^m h_j)$ is a tree.
The remaining $m$ generators of $\ub(G,2)$ correspond to cycles $\ti h_1,\dots,\ti h_m\subset G$
 passing through the selected edges $h_1,\dots,h_m$, respectively, when one robot stays 
 at a base point and the other robot moves along a cycle $\ti h_j$ without collisions.
\smallskip

\noindent
$\bullet$
Each relation says that motions of 2 robots along disjoint cycles commute.
\end{proposition}
\begin{proof}
By Subdivision Theorem~\ref{the:Subdivision} to compute the 2-point braid group $\ub(G,2)$,
 we may assume that $G$ has no loops and multiple edges 
 removing extra trivial vertices of degree~2.
Induction on the first Betti number $m$.
Base $m=0$ is Lemma~\ref{lem:2pointTreeGroupsUnordered}, 
 where every generator $\de_1\de_j^{-1}$ coming from 
 Proposition~\ref{pro:AddHangingEdgeUnordered} is represented 
 by a loop swapping 2 robots near a vertex of degree at least 3
 as shown in Figure~\ref{fig:PermuteRobotsTriod}.
\smallskip

In the induction step, for an edge $h\subset G$ from Lemma~\ref{lem:ChooseEdge},
 we show how a presentation of $\ub(G,2)$ differs from a presentation of 
 $\ub(G-h,2)$ satisfying the conditions by the induction hypothesis.
Since all cycles of $G-\nb(h)$ are in one connected component then $k=1$ 
 in Proposition~\ref{pro:CreateCyclesUnordered} and we skip the index $j$. 
So we add 1 new generator $\e(b\times h)\tau^{-1}$ that conjugates 
 the loops $\e(\bb\{x_n=r\})\e^{-1}$ and $\tau(\bb\{x_n=w\})\tau^{-1}$.
Geometrically, $\e(b\times h)\tau^{-1}$ represents 
 a motion when the 1st robot stays away from the 2nd robot 
 that completes a cycle $\ti h\subset G$ containing $h$.
\smallskip

It remains to show that the loops $(\bb\{x_n=r\})$ and $(\bb\{x_n=w\})$ 
 are homotopic, i.e. the new relator is a commutator.
Take the cycle $C\supset h$ from the construction of Lemma~\ref{lem:ChooseEdge}.
Since $C$ does not meet all cycles from $G-\nb(h)$, then we may move the 2nd robot
 along $C-h$ from $r$ to $w$ without collisions with the 1st robot moving
 along the cycles $\bb$ generating $\pi_1(G-\nb(h))$.
This gives a free homotopy from $(\bb\{x_n=r\})$ to
 $(\bb\{x_n=w\})=(b\times(C-h))(\bb\{x_n=r\})(b\times(C-h))^{-1}$.
\smallskip

During the motion $(b\times(C-h))$ the 1st robot is fixed at the base point 
 $b$ in $G-\nb(h)$, the 2nd moves along $C-h$ avoiding all cycles of $G-\nb(h)$.
In Proposition~\ref{pro:CreateCyclesUnordered} we may choose the path $\tau$ 
 from $a$ to $b\times w$ in $\UD(G-h,2)$ so that $\tau=\e\cdot(b\times(C-h))$. 
Then the loops $\e(\bb\{x_n=r\})\e^{-1}$ and $\tau(\bb\{x_n=w\})\tau^{-1}$
 are homotopic with the fixed base point $a\in\UD(G-h,2)$. 
\end{proof}


\subsection{Further open problems}
\label{subs:OpenProblems}
\noindent
\smallskip

Generalising the results of sections~\ref{sect:FundamentalGroupsUnordered}  
 and \ref{sect:ComputingGroupsUnordered} to ordered robots is left to followers.

\begin{problem}
Design and implement an algorithm computing a presentation of
 the pure braid group of an arbitrary connected graph similarly to 
 Algorithm~\ref{alg:GraphBraidGroups}.
\end{problem}

Our experience shows that presentations of planar graph braid groups may naturally 
 contain relators that are not commutators if there are no enough disjoint cycles.
So we state the problem opposite to \cite[Conjecture~5.7]{FS2} saying
 that all 2-point braid groups of planar graphs have presentations 
 where all relators are commutators. 

\begin{problem}
Check the conjecture that if $\ub(G,n)$ has a presentation such that 
 all relators are commutators then G can be constructed 
 as in Lemma~\ref{lem:ChooseEdge}.
\end{problem}


\end{document}